# Multi-dimensional filtering: Reducing the dimension through rotation

Julia Docampo Sánchez*¶    Jennifer K. Ryan*§¶

Mahsa Mirzargar†    Robert M. Kirby‡

October 7, 2016

## Abstract

Over the past few decades there has been a strong effort towards the development of Smoothness-Increasing Accuracy-Conserving (SIAC) filters [21] for Discontinuous Galerkin (DG) methods, designed to increase the smoothness and improve the convergence rate of the DG solution through this post-processor. These advantages can be exploited during flow visualization, for example by applying the SIAC filter to the DG data before streamline computations [Steffan et al., IEEE-TVCG 14(3): 680-692]. However, introducing these filters in engineering applications can be challenging since a tensor product filter grows in support size as the field dimension increases, becoming computationally expensive. As an alternative, [Walfisch et al., JOMP 38(2);164-184] proposed a univariate filter implemented along the streamline curves. Until now, this technique remained a numerical experiment. In this paper we introduce the SIAC line filter and explore how the orientation, structure and filter size affect the order of accuracy and global errors. We present theoretical error estimates showing how line filtering preserves the properties of traditional tensor product filtering, including smoothness and improvement in the convergence rate. Furthermore, numerical experiments are included, exhibiting how these filters achieve the same accuracy at significantly lower computational costs, becoming an attractive tool for the scientific visualization community.

*School of Mathematics, University of East Anglia, Norwich NR4 7TJ, United Kingdom
†Department of Computer Science, University of Miami, Coral Gables, FL 33124, USA
‡School of Computing, University of Utah, Salt Lake City, UT 84112, USA
§Corresponding author. Email: Jennifer.Ryan@uea.ac.uk

¶The authors are sponsored in part by the Air Force Office of Scientific Research (AFOSR), Computational Mathematics Program (Program Manager: Dr. Jean-Luc Cambier), under grant number FA8655-13-1-3017.

# 1 Introduction

This paper presents a new computationally efficient filtering technique employed to improve the quality of *multi-dimensional* numerical solutions obtained through Discontinuous Galerkin (DG) methods. We do this using a new approach to Smoothness-Increasing Accuracy-Conserving (SIAC) Filtering, which we call SIAC Line filtering. SIAC Filters [21] are a post-processing technique designed to accelerate the convergence rate and increase the smoothness of DG solutions. Traditional applications of SIAC filters require a tensor product construction. Here, we theoretically and computationally demonstrate that using an appropriate rotation of the one-dimensional SIAC filter, we can preserve the properties of the original post-processor. In fact, our computational results show that in addition to smoothness recovery, the new solution is generally more accurate than the original DG solution. Additionally, the theory provides guidelines on choosing the appropriate rotation, which is linked to the underlying divided-difference estimates based on the DG mesh.

Using a one-dimensional idea for multi-dimensional data makes SIAC an attractive tool for the visualization community. Flow visualization through particle tracking methods such as streamlines and streaklines is a common technique used to provide insight into fluid dynamics. Among the many techniques used for Computational Fluid Dynamics, DG methods are one family of numerical schemes that allow for generating data for flow visualization. They are robust, high order methods which can handle complicated geometries as well as effectively solve solutions containing shocks [6]. DG schemes, like Finite Element (FEM) and Finite Volume (FVM) Methods, use a variational form to solve Partial Differential Equations (PDEs). However, unlike FEM that require global continuity, a DG solution is continuous only inside the elements. The solution across the element interface is controlled through a numerical flux that is only weakly continuous; as a result, the error exhibits high frequency oscillations. Hence, visualizing DG solutions can be challenging. The numerical solution has low levels of continuity and most visualization techniques assume smooth field conditions. In order to increase the levels of continuity, typically a post-processor is implemented.

There has been ongoing work on the application of SIAC filters for DG solutions to improve the flow conditions where streamlines are subsequently computed. The authors of [7] implemented the traditional multidimensional filter obtained as a tensor product of univariate filters along each Cartesian axis. This configuration, a natural extension of the one dimensional case, allows for proving error estimates both for uniform and nonuniform cases [5]. However, the foundations for proving superconvergence assume only smooth initial data and link the filter directly to the underlying mesh, restricting the choices on the area of the domain from which information is extracted. We hypothesized that changing the direction in which information was filtered could improve the results. Therefore, we considered rotated filters: SIAC filters that were no longer Cartesian coordinate aligned and had variable orientation. This idea comes from a visualization perspective, questioning if orienting the filter with the flow direction and changing the support size plays a role in improving the quality of the filtered solution.

Making a filter viable for flow visualization comes with computational challenges because it requires robustness, relatively low computational intensity, and short simulation times. In [27], a numerical experiment on streamline visualization was performed. In order to save computational costs and avoid using a 2D filter, the authors applied a type of one dimensional filter along the streamline curves. However, the theoretical and numerical investigation into the effectiveness of these filters on typical test problems was not carried out. Here we perform this investigation on *SIAC Line Filters* and show that this new approach is a computationally efficient technique for post-processing multidimensional fields that uses only one dimension. This family of filters transforms the 2D integral of the convolution into a line integral. Hence, from a computational point of view, the advantages are immediate. Furthermore, we prove that it is possible to extract superconvergence for such filters and present numerical results supporting the theory. The results of the line filter are compared to the original 2D Cartesian coordinate aligned SIAC filter.

The results in this paper combine previous investigations in SIAC Filtering. SIAC filters have traditionally been used to reduce the error oscillations and recover smoothness in the solution and its derivatives [14, 12, 15, 24, 22, 20, 11]. The filters were originally designed for accuracy enhancement of FEMs [2, 17] and later applied to DG [5]. The post-processor extracts the hidden "superconvergence" of these methods; for linear hyperbolic problems, the filtered solution is of order $2k+1$, where $k$ denotes the polynomial space degree used for the DG approximation which is order $k+1$ convergent. Hence, in addition to increasing the smoothness, for smooth initial data and linear problems, the filtered solution is generally more accurate than the DG solution.

We will begin by briefly reviewing the DG method and the original post-processor. In Section 3 we provide details on how to rotate the filter and reduce the dimension for multi-dimensional filtering. We additionally provide theoretical error estimates for linear hyperbolic equations. In Section 4 numerical experiments showing smoothness recovery and accuracy enhancement are provided and a discussion of the computational benefits of this type of post-processing are provided in Section 5. Our conclusions are formulated in Section 6 to close this article.

## 2 Background

The theory of SIAC filtering for DG methods relies on the divided differences of the numerical solution. Using a piecewise polynomial basis of degree $k$, the numerical solution is typically of order $k+1$ under the $L^2$ norm in both the approximation and divided differences for linear hyperbolic equations. However, DG solutions have "hidden" superconvergence. In [5] it was proven that the DG approximation has $2k+1$ convergence in the negative-order norm for the approximation and the divided differences. SIAC filters exploit this fact and can achieve $2k+1$ order in the $L^2$ norm for the actual solution. In order to understand how we can extract superconvergence, we will begin by introducing the DG scheme and the theoretical error estimates.

### 2.1 The DG Scheme and its Divided Differences

The first DG method was developed to solve the neutron transport equation in 1973 by Reed an Hill [19]. Today, these methods extend to many types of PDEs.

Examples are solving non-linear combined problems such as the incompressible and compressible Navier-Stokes equations [1, 9].

In this paper, we concentrate on linear hyperbolic conservation laws and use the advection equation as the model problem. DG schemes for such problems have been studied in depth by [4, 6] and here we only address the basic idea of the scheme. Consider the linear hyperbolic problem:

$$\begin{cases} u_t + \sum_{i=1}^{d} A_i u_{x_i} + A_0 u = 0, \quad (\mathbf{x}, t) \in \Omega \times [0, T] \\ u(\mathbf{x}, 0) = u_0, \end{cases} \tag{1}$$

where $A_i$, $i = 1 \ldots, d$ are linear, $\mathbf{x} = (x_1, x_2, \ldots, x_d)$ and $u$ represents the advection of the conserved quantity.

The first step of the DG method is to choose a suitable tessellation $\mathcal{T}(\Omega) = \sum e$ of the domain $\Omega$ and a piece-wise polynomial approximation space:

$$V_h^k = \left\{ v \in L^2(\Omega) : v \in \mathbb{P}^k(e), \ \forall e \in \mathcal{T}(\Omega) \right\}.$$

Then, the DG solution is obtained using the variational form of Equation (1). It is the unique function $u_h \in V_h^k$ satisfying

$$\int_e (u_h)_t v dx - \sum_{i=1}^{d} \left( \int_e A_i u_h(x,t) v_{x_i} dx \right) + \int_e A_0 u_h v dx + \sum_{i=1}^{d} \int_{\partial e} \widehat{A_i u_h} \cdot \mathbf{n} v dS = 0 \tag{2}$$

for all $v \in V_h^k$ and for every element of the tessellation. The term $\widehat{A_i u_h}$ refers to the numerical flux, the function enforcing weak continuity across the element interfaces, which is typically taken to be the upwind flux.

**Theorem 2.1** ([5]). *Let $u$ be the exact solution and $u_h$ the DG approximation to the Boundary Value Problem (1) with periodic boundary conditions. For a uniform mesh, we obtain the following error estimates:*

$$\|\partial_h^\alpha (u - u_h)\|_{0,\Omega} \leq Ch^{k+1} \tag{3}$$

*in the $L^2$-norm and in the negative order norm:*

$$\|\partial_h^\alpha (u - u_h)\|_{-(k+1),\Omega} \leq Ch^{2k+1}, \tag{4}$$

where

$$\|u\|_{-\ell,\Omega} = \sup_{\phi \in \mathcal{C}_0^\infty(\Omega)} \frac{(u,\phi)_\Omega}{\|\phi\|_{\ell,\Omega}}, \quad \|\phi\|_{\ell,\Omega} = \left( \sum_{|\alpha| \leq \ell} \|D^\alpha u\|_\Omega^2 \right)^{\frac{1}{2}} \quad \text{and } \ell > 0.$$

Here $k$ denotes the polynomial order used for the DG approximation, $\alpha = (\alpha_1, \alpha_2, \ldots, \alpha_d)$ a multi-index. It is also necessary to introduce $\partial_h^\alpha$, the (scaled) divided difference:

$$\partial_h^\alpha = \partial_{h,1}^{\alpha_1} \partial_{h,2}^{\alpha_2} \cdots \partial_{h,d}^{\alpha_d} \quad \partial_{h,j} f(x) = \frac{1}{h}(f(x + h/2) - f(x - h/2)), \tag{5}$$

$$\partial_{h,j}^{\alpha_j} f = \partial_{h,j}(\partial_{h,j}^{\alpha_j - 1} f), \quad \alpha_j > 1, j = 1, \ldots, d. \tag{6}$$

For the purposes of this article, we take $d = 2$.

Finally, we introduce the Lemma that allows us to switch between the $L^2$ and the negative-order norms.

**Lemma 2.1** *(Bramble and Schatz [2]).* *Let $\Omega_0 \subset\subset \Omega_1 \subset\subset \Omega$, $\Omega$ bounded domain in $\mathbb{R}^d$ and $s$ be an arbitrary but fixed nonnegative integer. Then, for $u \in H^s(\Omega_1)$, there is a constant $C$ such that*

$$||u||_{0,\Omega_0} \leq C \sum_{|\alpha| \leq s} ||D^\alpha u||_{-s,\Omega_1}. \tag{7}$$

In the next section, we will see that the B-Spline kernel transforms the differential operator $D^\alpha$ into a divided differences operator. This allows the use of the negative order estimate of Theorem 2.1 over the filtered solution, giving $2k+1$ accuracy in the $L^2$ norm [5].

## 2.2 SIAC Filters

Before introducing the filter rotation, we briefly review the original post-processor from which it derives. For a much more detailed description on the properties and implementation of SIAC filters, we refer the reader to [13, 14, 21, 16].

The post-processor is a continuous convolution:

$$u_h^\star(x,T) = \int_{-\infty}^{\infty} K_H^{(2k+1,k+1)}(x-y) u_h(y,T)\, dy, \quad x \in \Omega \tag{8}$$

where $u_h$ denotes the DG solution at final time and the kernel is a linear combination of central B-Splines:

$$K^{(2k+1,k+1)}(\eta) = \sum_{\gamma=-k}^{k} c_\gamma \psi^{(k+1)}(\eta - \gamma). \tag{9}$$

Here, $\gamma$ denotes the B-Splines centres. The kernel subindex $H$ in equation (8) acts as a scaling factor, *i.e.*, $K_H(x-y) = \frac{1}{H} K\left(\frac{x-y}{H}\right)$. To give an idea of the filter size, for uniform meshes, the usual scaling choice is $H = h$, where $h$ denotes the element size. The superindexes $(2k+1, k+1)$ indicate the number of B-Splines used to build the kernel $(2k+1)$ and the spline order $(k+1)$. Basis Splines (B-Splines) are local functions providing maximum approximation order with minimum support. The central B-Splines are a particular case which are computationally very atractive since they can be calculated using the recurrence formula

$$\psi^{(1)}(x) = \chi_{[-1/2, 1/2]}(x) \tag{10}$$

$$\psi^{(k+1)}(x) = \psi^{(k)} \star \psi^{(1)}(x) \tag{11}$$

$$= \frac{1}{k}\left(\left(\frac{k+1}{2} + x\right)\psi^{(k)}\left(x + \frac{1}{2}\right) + \left(\frac{k+1}{2} - x\right)\psi^{(k)}\left(x - \frac{1}{2}\right)\right). \tag{12}$$

Moreover, these splines have the following property for the derivatives:

$$d^\alpha \psi^{(k+1)} = \partial_{h=1}^\alpha \psi^{(k+1-\alpha)}, \quad \partial_h^\alpha = \alpha^{th} \text{ divided difference.} \tag{13}$$

We will not give further details on these spline functions and suggest [8] and [25] for a complete description. Finally, the kernel coefficients, $c_\gamma$, dictate each of the B-Spline weights and are determined by imposing the polynomial reproduction property

$$K^{(2k+1,k+1)} \star x^p = x^p, \ p = 0, \ldots, 2k. \tag{14}$$

The kernel presented in equation (9) is symmetric in the sense that the support is centred around the post-processing point and it expands equally in every direction. There are alternative kernel versions, giving rise to one-sided [23] and position dependent [26, 11] SIAC filters, and more recently, the non-uniform knot based PSIAC filters [18]. These kernels include a shifting parameter in the B-Splines, translating the support towards one direction. In Figure 1 we show the B-Spline functions together with a symmetric and one-sided kernel to ilustrate the difference. These alternative versions attempt to address issues related to domain boundaries and near-shock regions. Since beyond the computational domain there is no information, near the boundaries the symmetric kernel can not be implemented. Instead, it is replaced by a boundary filter, allowing post-processing points by pushing the support towards the interior of the domain. Furthermore, for solutions containing shocks, taking information near the shock may produce an undesirable smooth region. However, in this paper we will not tackle these problems and concentrate only on the symmetric filter, assuming periodic conditions and linear hyperbolic problems.

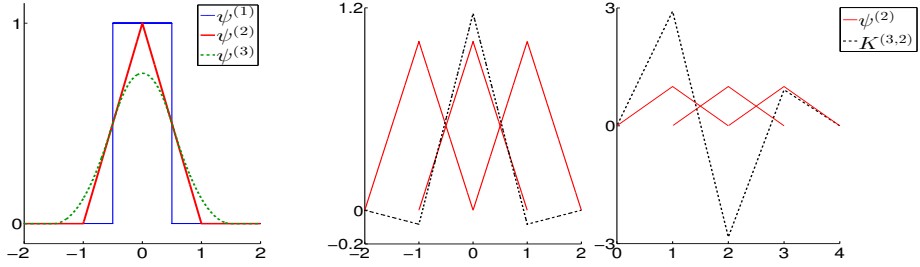

Figure 1: B-Splines (left) and the symmetric (centre) and one-sided RS [23] (right) kernels.

**Theorem 2.2** *(Cockburn, Luskin, Shu, and Süli [5]) Under the same conditions in Theorem 2.1 and if*

$\Omega_0 + 2supp\left(K_h^{(2k+1,k+1)}\right) \subset\subset \Omega_1 \subset \Omega$, *then for $H = h$ (h mesh size):*

$$\left\| u - K_h^{(2k+1,k+1)} \star u_h \right\|_{0,\Omega_0} \leq Ch^{2k+1}. \tag{15}$$

Here we sketch the proof in order to illustrate the important components for proving the same properties for the rotated filter.

$$\left\| u - K_h^{(2k+1,k+1)} \star u_h \right\|_{0,\Omega_0} \leq \underbrace{\left\| u - K_h^{(2k+1,k+1)} \star u \right\|_{0,\Omega_0}}_{\Theta_1} + \underbrace{\left\| K_h^{(2k+1,k+1)} \star (u - u_h) \right\|_{0,\Omega_0}}_{\Theta_2}.$$

The term $\Theta_1$ is bounded using property (14), polynomial reproduction. The second term relies on property (13), the ability to switch the derivative to a divided difference. For multiple dimensions, this relies on $\alpha = (\alpha_1, \alpha_2, \ldots, \alpha_d)$ and therefore contains a multi-dimensional derivative.

□

**Remark 2.1** *The polynomial reproduction property implies that convolving the exact solution with the filter produces an error of order $\mathcal{O}(h^{2k+1})$, with $2k$ being maximum polynomial degree of reproduction. This is controlled by the number of B-Splines used during kernel construction.*

**Remark 2.2** *The divided differences play a key role for bounding the error component corresponding to the filtered DG approximation. The $2k+1$ accuracy is achieved by virtue of Theorem 2.1 using the B-Spline derivative property given in equation* (13).

## 3 SIAC Line Filters

Before introducing Smoothness-Increasing Accuracy-Conserving Line Filters, we introduce the coordinate rotation necessary for these filters in the context of the two-dimensional Cartesian axis aligned filter. We then show how the theory automatically translates to this one-dimensional kernel for multi-dimensional data.

### 3.1 Rotating the Kernel support

Now we discuss rotating the two-dimensional SIAC kernel. The post-processor aligned with the Cartesian axis is built as a tensor product of univariate kernels:

$$u^\star(\overline{x},\overline{y}) = \int_{-\infty}^{\infty} \int_{-\infty}^{\infty} K_{H_x}^{(2k+1,k+1)}(\overline{x}-x) K_{H_y}^{(2k+1,k+1)}(\overline{y}-y) u_h(x,y)\, dx\, dy. \quad (16)$$

For uniform meshes, provided the kernel scaling is of the form $H = mh$, where $m \in \mathbb{Z}^+$ and $h$ is the mesh size, it is possible to show $2k+1$ accuracy. However, as soon as the mesh uniformity assumption drops, finding a suitable scaling becomes complicated. A detailed theoretical discussion on the kernel scaling and nonuniform meshes can be found in [10, Ch. 4] and see [7] for a numerical study.

The idea of allowing kernel rotations that change the support orientation comes from practical applications of SIAC filters. In terms of robustness, formulation (16) is restrictive in the sense that there is only one possible choice for the kernel support: a box aligned with the cartesian axis. If one wants to implement the post-processor during flow visualization, it is reasonable to ask whether keeping the kernel aligned with the mesh is more relevant than aligning the kernel with the flow direction. This argument does not contradict any of the previous work on SIAC filtering, since the aim of post-processing was extracting superconvergence. We wish to highlight that *ensuring superconvergence does not necessarily imply ensuring error minimization*. This fact can be observed in the numerical results presented in Section 4.

The rotated filter consists of rewriting the convolution in a new basis $\mathcal{B}_2 := \left\{\vec{k_x}, \vec{k_y}\right\}$, given by the rotation angle and the kernel direction vectors: $\vec{k_x} = (\cos\theta,\ \sin\theta)$, $\vec{k_y} = (\cos\theta + \pi/2,\ \sin\theta + \pi/2) = (-\sin\theta, \cos\theta)$. Identifying $\mathcal{B}_1 := \{\mathbf{e_1} = (1,0), \mathbf{e_2} = (0,1)\}$ with the cartesian basis and using the **change-of-basis** matrices:

$$P_{\mathcal{B}_2 \leftarrow \mathcal{B}_1} = \begin{pmatrix} \cos\theta & -\sin\theta \\ \sin\theta & \cos\theta \end{pmatrix} \quad \text{and} \quad P_{\mathcal{B}_1 \leftarrow \mathcal{B}_2} = \begin{pmatrix} \cos\theta & \sin\theta \\ -\sin\theta & \cos\theta \end{pmatrix} \quad (17)$$

we can exchange the coordinates of a point from each reference system:

$$X = (x_1, x_2)_{\mathcal{B}_1} = P_{\mathcal{B}_1 \leftarrow \mathcal{B}_2} \cdot X' \Leftrightarrow X' = (x'_1, x'_2)_{\mathcal{B}_2} = P_{\mathcal{B}_2 \leftarrow \mathcal{B}_1} \cdot X.$$

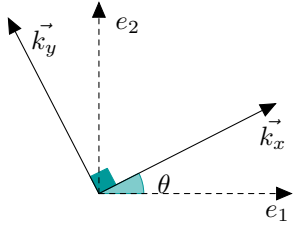

Figure 2: Illustration of the two coordinate systems basis vectors.

The filtering convolution is defined in the new coordinate system:

$$u^\star(\overline{x}, \overline{y}) = \int_{-\infty}^{\infty} \int_{-\infty}^{\infty} K_{H_x}^{(2k+1,k+1)}(\overline{x} - x') K_{H_y}^{(2k+1,k+1)}(\overline{y} - y') u_h(x', y') dx' dy' \quad (18)$$

and each (symmetric) kernel is defined using (9).

**Remark 3.1** *This definition is consistent with the original post-processor. This can be seen using zero rotation. It is the particular case when $x = x'$ $(y = y')$.*

### 3.2 Reducing the dimension through line filtering

In Theorem 2.2 we highlighted the important role played by the divided differences for proving superconvergence of the filtered solution; the proof relies on the ability of the kernel to transfer the derivatives to the DG approximation as divided differences and then apply Theorem 2.1. With an axis aligned kernel, a tensor product construction,

$$K_H^{(2k+1,k+1)}(x,y) = K_{H_x}^{(2k+1,k+1)}(x) \otimes K_{H_y}^{(2k+1,k+1)}(y),$$

is necessary in order to compute the multi-dimensional derivatives:

$$D^\alpha K_H^{(2k+1,k+1)}(x,y) = \frac{d^{\alpha_1}}{dx} K_{H_x}^{(2k+1,k+1)}(x) \frac{d^{\alpha_2}}{dy} K_{H_y}^{(2k+1,k+1)}(y), \quad \alpha_1 + \alpha_2 = \alpha.$$

On the other hand, rotating the kernel produces a great advantage since a single kernel direction allows for differentiation in terms of the original basis under all variables. For example, consider the rotated kernel direction $k_x$ in equation (18) and assume that $\theta \in (0, \pi/2)$:

$$K_{H_x}^{(2k+1,k+1)}(x') = K_{H_x}^{(2k+1,k+1)}(x,y), \quad \text{since } x' = (\cos\theta \ -\sin\theta) \begin{pmatrix} x \\ y \end{pmatrix}. \quad (19)$$

Exploiting this fact we can avoid tensor products and reduce the filter dimension, transforming the convolution into a line integral whilst preserving the 2D SIAC properties. Therefore, we only need to show that the kernel derivatives can still be expressed as a combination of divided differences in the $x$ and $y$ directions and then the same error estimates will hold.

We would like to emphasize that the idea of univariate SIAC filters for multi-dimensional domains was first introduced in [27]. They showed the potential of this technique with an empirical study on streamlines, implementing a one-sided filter along the curve using arc-length parametrization. Here we mathematically develop that idea and define **SIAC line filters**: the family of rotated SIAC filters with support expanding only along a segment inside the 2D domain. We will define the line kernel together with the B-Splines and formalise this lower dimension filtering approach presenting theoretical error estimates. We conclude this section with supporting numerical results showing both smoothness recovery and superconvergence.

**Definition 3.1** (*SIAC line kernels*). *Let $\Gamma \subset \mathbb{R}^2$ be the line parametrized by the arc length*
$$\Gamma(t) = t(\cos\theta, \sin\theta) \qquad t \in \mathbb{R},\ \theta\ \text{fixed,} \tag{20}$$
*and its inverse*
$$\Gamma^{-1}(x,y) = x\cos\theta + y\sin\theta. \tag{21}$$
*Then the B-Spline along the $\Gamma$ line is defined by:*
$$\tilde{\psi}_\theta^{(\ell)}(x,y) = \begin{cases} \psi^{(\ell)}\left(\Gamma^{-1}(x,y)\right) & \text{if } (x,y) \in \Gamma(t) \\ 0 & \text{otherwise,} \end{cases} \tag{22}$$
*and has compact support*
$$\text{supp}\ \tilde{\psi}_\theta^{(\ell)} = (t\cos\theta, t\sin\theta), \quad t \in \left[-\frac{\ell}{2}, \frac{\ell}{2}\right]. \tag{23}$$

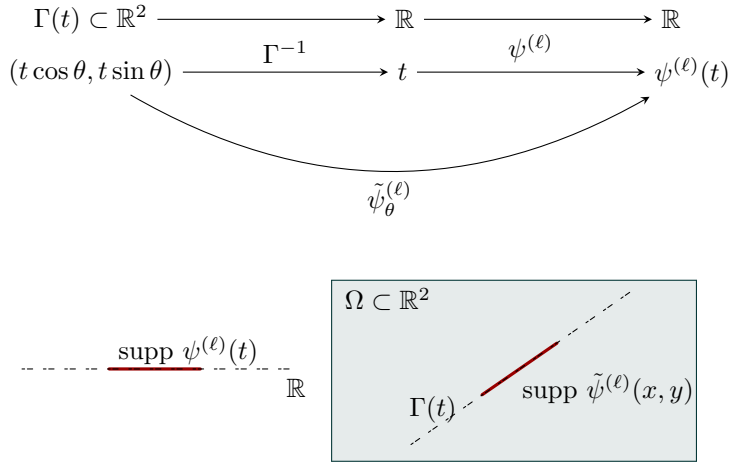

Figure 3: Illustration of an univariate B-Spline support along a line in $\mathbb{R}^2$.

The **SIAC line kernel** is constructed as a linear combination of these (scaled) B-Splines and the symmetric version, originally introduced in [2, 5], has the following formula:

$$K_{H,\Gamma}^{(2k+1,k+1)}(t) = \sum_{\gamma=-k}^{k} c_\gamma \psi_{\theta,H}^{(k+1)}(t-\gamma) \qquad (24)$$

in arc length coordinates, or alternatively by:

$$K_{H,\Gamma}^{(2k+1,k+1)}(x,y) = \sum_{\gamma=-k}^{k} c_\gamma \tilde{\psi}_{\theta,H}^{(k+1)}\left(\Gamma^{-1}\left(x-\gamma\cos\theta, y-\gamma\sin\theta\right)\right) \qquad (25)$$

in the Cartesian system.

In order to characterize the derivatives of such B-Splines we introduce a particular type of divided differences.

**Definition 3.2** (Directional Divided Difference). Consider the direction given by the vector $\vec{u} = (u_x, u_y)$. Then the scaled directional divided difference with respect to $\vec{u}$ is defined by

$$\partial_{\vec{u},H} f(x,y) = \frac{1}{H}\left(f\left(x+\frac{H}{2}u_x, y+\frac{H}{2}u_y\right) - f\left(x-\frac{H}{2}u_x, y-\frac{H}{2}u_y\right)\right). \qquad (26)$$

The $\alpha$-directional divided difference is then given by

$$\partial_{\vec{u},H}^\alpha f(x,y) = \partial_{\vec{u},H}\left(\partial_{\vec{u},H}^{\alpha-1} f(x,y)\right), \quad \alpha > 1. \qquad (27)$$

**Lemma 1** For a B-Spline satisfying Definition 3.1, the (scaled) directional divided differences along the line $\Gamma$ are equal to the (scaled) divided differences of the univariate B-Spline along the arc parameter, ie,

$$\partial_{u_\theta,H}^\alpha \tilde{\psi}_{\theta,H}^{(\ell-\alpha)}(x,y) = \partial_H^\alpha \psi_H^{(\ell-\alpha)}(t), \quad u_\theta = (\cos\theta, \sin\theta). \qquad (28)$$

*Proof.* Consider $\alpha = 1$:

$$\partial_{u_\theta,H} \tilde{\psi}_{\theta,H}^{(\ell-1)}(x,y) = \frac{\tilde{\psi}_{\theta,H}^{(\ell-1)}\left(x+\frac{H}{2}\cos\theta, y+\frac{H}{2}\sin\theta\right)}{H} - \frac{\tilde{\psi}_{\theta,H}^{(\ell-1)}\left(x-\frac{H}{2}\cos\theta, y-\frac{H}{2}\sin\theta\right)}{H}$$

Examining the first term:

$$\frac{\tilde{\psi}_{\theta,H}^{(\ell-1)}\left(x+\frac{H}{2}\cos\theta, y+\frac{H}{2}\sin\theta\right)}{H} = \frac{\psi_H^{(\ell-1)}\left(\Gamma^{-1}\left(x+\frac{H}{2}\cos\theta, y+\frac{H}{2}\sin\theta\right)\right)}{H}$$
$$= \frac{\psi_H^{(\ell-1)}\left(t+\frac{H}{2}\right)}{H}.$$

Hence,

$$\partial_{u_\theta,H} \tilde{\psi}_{\theta,H}^{(\ell-1)}(x,y) = \frac{\psi_H^{(\ell-1)}\left(t+\frac{H}{2}\right)}{H} - \frac{\psi_H^{(\ell-1)}\left(t-\frac{H}{2}\right)}{H} = \partial_H \psi^{(\ell)}(t).$$

Using induction, the higher order divided differences in equation (28) can be obtained. □

**Lemma 2** *The $\alpha$−derivative of the B-Spline from Definition 3.1 can be expressed as a sum of $\alpha$-directional divided differences using the basis vectors:*

$$u_x^\theta = (\cos\theta, 0), \quad u_y^\theta = (0, \sin\theta)$$

*through the formula*

$$D^\alpha \tilde{\psi}_{\theta,H}^{(\ell)}(x,y) = \sin^{\alpha_1}\theta \cos^{\alpha_2}\theta \sum_{m=0}^{\alpha} \binom{\alpha}{m} \partial_{u_x^\theta,H}^{\alpha-m} \partial_{u_y^\theta,H}^{m} \tilde{\psi}_{\theta,H}^{(\ell-\alpha)}\left(x - \frac{m}{2}H\cos\theta, y + \frac{\alpha-m}{2}H\sin\theta\right)$$

*where $\alpha_1 + \alpha_2 = \alpha$.*

*Proof.* Differentiating a B-Spline of degree $\ell$ gives

$$D^\alpha \tilde{\psi}_{\theta,H}^{(\ell)}(x,y) = \sin^{\alpha_1}\theta \cos^{\alpha_2}\theta \left(\frac{d^\alpha \psi^{(\ell)}(t)}{dt^\alpha}\right) = \sin^{\alpha_1}\theta \cos^{\alpha_2}\theta \left(\partial_H^\alpha \psi_H^{(\ell-\alpha)}(t)\right), \quad \alpha_1 + \alpha_2 = \alpha.$$

Then use Lemma 1 to obtain:

$$D^\alpha \tilde{\psi}_{\theta,H}^{(\ell)}(x,y) = \sin^{\alpha_1}\theta \cos^{\alpha_2}\theta \left(\partial_{u_\theta,H}^\alpha \tilde{\psi}_{\theta,H}^{(\ell-\alpha)}(x,y)\right).$$

In order to show that

$$\partial_{u_\theta,H}^\alpha \tilde{\psi}_{\theta,H}^{(\ell-\alpha)}(x,y) = \sum_{m=0}^{\alpha} \binom{\alpha}{m} \partial_{u_x^\theta,H}^{\alpha-m} \partial_{u_y^\theta,H}^{m} \tilde{\psi}_{\theta,H}^{(\ell-\alpha)}\left(x - \frac{m}{2}H\cos\theta, y + \frac{\alpha-m}{2}H\sin\theta\right),$$

we use induction. Consider $\alpha = 1$:

$$\partial_{u_\theta,H} \tilde{\psi}_{\theta,H}^{(\ell-1)}(x,y) = \frac{1}{H}\left(\tilde{\psi}_{\theta,H}^{(\ell-1)}\left(x + \frac{H}{2}\cos\theta, y + \frac{H}{2}\sin\theta\right) - \tilde{\psi}_{\theta,H}^{(\ell-1)}\left(x - \frac{H}{2}\cos\theta, y - \frac{H}{2}\sin\theta\right)\right).$$

Adding and subtracting the term

$$\frac{1}{H}\tilde{\psi}_{\theta,H}^{(\ell-1)}\left(x - \frac{H}{2}\cos\theta, y + \frac{H}{2}\sin\theta\right)$$

gives:

$$\partial_{u_\theta,H} \tilde{\psi}_{\theta,H}^{(\ell-1)}(x,y) = \partial_{u_x^\theta,H} \tilde{\psi}_{\theta,H}^{(\ell-1)}\left(x, y + \frac{H}{2}\sin\theta\right) + \partial_{u_y^\theta,H} \tilde{\psi}_{\theta,H}^{(\ell-1)}\left(x - \frac{H}{2}\cos\theta, y\right).$$

For the general case, we can write:

$$\partial_{u_\theta,H}^\alpha \tilde{\psi}_{\theta,H}^{(\ell-\alpha)}(x,y) = \partial_{u_\theta,H}\left(\partial_{u_\theta,H}^{\alpha-1} \tilde{\psi}_{\theta,H}^{(\ell-(\alpha-1))}(x,y)\right)$$

and the formula can be shown by induction using a similar proof to the binomial theorem. The full proof is shown in the appendix. □

**Remark 3.2** *Lemma 2 shows that the B-Spline derivatives can be computed as a sum of directional divided differences along the standard basis vectors.*

We now introduce our main Theorem that allows us to give error estimates for SIAC line filters.

**Theorem 3.1** *Let $u$ be the exact solution to Problem (1) with $d = 2$ and periodic boundary conditions. Let $u_h$ be the DG approximation over a uniform mesh and denote by $h_x$ and $h_y$ the mesh size. Consider the line filter $K_{\Gamma,H}^{(2k+1,k+1)}$ along $\Gamma(t) = t(\cos\theta, \sin\theta)$, $\theta$ fixed. If $\theta = \arctan\left(\frac{h_y}{h_x}\right)$ and $H = h_x \cos\theta + h_y \sin\theta$, then:*

$$\left\| u - K_{\Gamma,H}^{(2k+1,k+1)} \star u_h \right\|_{0,\Omega_0} \leq C h^{2k+1}. \tag{29}$$

*Proof.*

$$\left\| u - K_{\Gamma,H}^{(2k+1,k+1)} \star u_h \right\|_{0,\Omega_0} \leq \underbrace{\left\| u - K_{\Gamma,H}^{(2k+1,k+1)} \star u \right\|_{0,\Omega_0}}_{\Theta_1} + \underbrace{\left\| K_h^{(2k+1,k+1)} \star (u - u_h) \right\|_{0,\Omega_0}}_{\Theta_2}.$$

The line filter preserves the polynomial reproduction property:

$$K_{\Gamma,H}^{(2k+1,k+1)} \star x^p = x^p, \quad p = 0, \ldots, 2k, \tag{30}$$

and therefore the first term is bounded similar to the first term in Theorem 2.2. For the second term, we need to show that the directional divided differences allow us to write an expression similar to

$$D^\alpha \left( K_H^{(2k+1,k+1)} \star (u - u_h) \right) = K_H^{(2k+1,k+1-\alpha)} \star \partial_H^\alpha (u - u_h)$$

in order to obtain a bound of the form of:

$$\Theta_2 \leq C_1 C_2 \sum_{|\alpha| \leq k+1} \|\partial_h^\alpha (u - u_h)\|_{-k+1,\Omega_1}. \tag{31}$$

Let $\ell = k + 1$ and denote the error by $e = u - u_h$. Notice that since the kernel is a linear combination of B-Splines, it is sufficient to study one B-Spline alone. Lemma 2 allows us to write

$$D^\alpha \tilde{\psi}_{\theta,H}^{(\ell)} \star e = \cos^{\alpha_1}\theta \sin^{\alpha_2}\theta \cdot \tilde{\psi}_{\theta,H}^{(\ell-\alpha)} \star \partial_{u_\theta,H}^\alpha e, \quad \alpha_1 + \alpha_2 = \alpha$$

$$= \sin^{\alpha_1}\theta \cos^{\alpha_2}\theta \tilde{\psi}_{\theta,H}^{(\ell)} \star \left( \sum_{m=0}^{\alpha} \binom{\alpha}{m} \partial_{u_x^\theta,H}^{\alpha-m} \partial_{u_y^\theta,H}^m e \left( x - \frac{m}{2} H \cos\theta, y + \frac{\alpha-m}{2} H \sin\theta \right) \right),$$

where $u_x^\theta = (\cos\theta, 0)$ and $u_y^\theta = (0, \sin\theta)$. The pair

$$\left\{ H = h_x \cos\theta + h_y \sin\theta, \ \theta = \arctan\left(\frac{h_y}{h_x}\right) \right\}$$

implies that the kernel scaling is also equivalent to

$$H = \frac{h_x}{\cos\theta} \quad \text{and} \quad H = \frac{h_y}{\sin\theta}.$$

This allows us to write the directional divided differences of the error function in the canonical basis $\mathcal{B}_1 = \{\mathbf{e}_1, \mathbf{e}_2\}$ using the mesh size:

$$\partial_{u_\theta^x,H} f(x,y) = \frac{1}{H}\left(f\left(x + \frac{H}{2}\cos\theta, y\right) - f\left(x - \frac{H}{2}\cos\theta, y\right)\right) = \frac{1}{H}\left(f\left(x + \frac{h_x}{2}\right) - f\left(x - \frac{h_x}{2}, y\right)\right)$$

$$= \cos\theta \cdot \partial_{\mathbf{e}_1,h_x} f(x,y).$$

Analogously: $\partial_{u^y_{\theta,H}} f(x,y) = \sin\theta \cdot \partial_{\mathbf{e}_2,h_y} f(x,y)$. Hence

$$D^\alpha \tilde{\psi}^{(\ell)}_{\theta,H} \star e = \sin^{\alpha_1+1}\theta \cos^{\alpha_2+1}\theta \tilde{\psi}^{(\ell)}_{\theta,H} \star \left( \sum_{m=0}^{\alpha} \binom{\alpha}{m} \partial^{\alpha-m}_{\mathbf{e}_1,h_x} \partial^m_{\mathbf{e}_2,h_y} e\left(x - \frac{m}{2}h_x, y + \frac{\alpha-m}{2}h_y\right) \right),$$

giving:

$$\Theta_2 \leq C_1 C_2 C(\theta) \sum_{|\alpha| \leq k+1} \left\| \sum_{m=0}^{\alpha} \binom{\alpha}{m} \partial^{\alpha-m}_{\mathbf{e}_1,h_x} \partial^m_{\mathbf{e}_2,h_y} e \right\|_{-k+1,\Omega_1}. \qquad (32)$$

The rest of the proof follows from Theorem 2.2. □

**Remark 3.3** *When a B-Spline is differentiated, as a consequence of the chain rule a $\sin\theta$ or $\cos\theta$ term appears. As a result, the constant term in equation (32) now includes the multiplying factor*

$$\sin^{\alpha_1+1}\theta \cos^{\alpha_2+1}$$

*which is always less than one (and decreasing with every power) since the rotation angle is defined by $\arctan(h_y/h_x)$. This means that the constant in front of equation (29) can actually be reduced. In the numerical experiments presented in the following section, we show cases where the line filter outperforms the 2D axis aligned filter.*

# 4 Numerical Results

The numerical experiments were done for the 2D advection equation:

$$\begin{cases} u_t + u_x + u_y = 0, & (x,y) \in [0, 2\pi]^2,\ t \in [0, T] \\ u_0(x,y) = u(x,0) \end{cases} \qquad (33)$$

with final time $T = 2$. The unfiltered solution was obtained using a DG scheme with an upwind flux over an uniform mesh. Two initial conditions were chosen:

1. $u_0(x,y) = \sin(x+y)$.
2. $u_0(x,y) = \sin x \cdot \cos y$.

In the following, we discuss the various aspects of filtering: smoothness and accuracy enhancement, including error reduction.

## 4.1 Recovering Smoothness

Here we discuss the potential of the line filter to recover smoothness. Since we are post-processing along a single direction, one can expect that it is only along that line where the filtered solution gains smoothness. The plots in Figures 4 and 5 show different error profiles corresponding to horizontal, vertical and diagonal cuts of the domain. We remind the reader that the kernel scaling $H$ is calculated through the formula $H = \mu h$, where $h$ is the DG mesh size. Notice that the $\theta = 3\pi/4$ and $\theta = \pi/4$ rotations use the value $\mu = \sqrt{2}$. This corresponds to the theoretical scaling $H = h(\cos(\theta) + \sin(\theta)) = \sqrt{2}h$. We observe that for the axis aligned kernel ($\theta = 0$), the filter gains smoothness along the filtering direction only. On the other hand, allowing a $\pi/4$ or $3\pi/4$ rotation produces a smooth profile in all three directions. This is because the rotated filter includes both the $x-$, $y-$ as well as mixed derivatives necessary for the theory to remain valid.

## 4.2 Accuracy enhancement

We now study the global error and convergence of the line filters. This study compares the line filtered solution with the DG solution. Whenever possible, we also include results obtained from filtering with the original tensor product filter aligned with the Cartesian axis. Due to a lack of computational resources, with the available machines we were forced to stop at $40 \times 40$ elements. However, since the numerical studies include three polynomial degrees, this suffices to give insight into the filter behaviour relative to each other. Table 1 shows the global $L^2$ errors and convergence rates for the 2D advection equation with initial condition $u_0(x,y) = \sin(x+y)$ using two different orientations of the SIAC line filter ($\theta = \pi/4$ and $\theta = 3\pi/4$) as well as the traditional axis aligned 2D filter. The $3\pi/4$ reduces the errors significantly over the unfiltered solution, even compared to the results from post-processing using a 2D Cartesian coordinate aligned kernel. Observing the order of accuracy, we see that the $\pi/4$ orientation has a faster convergence rate than the unfiltered solution and the errors are reduced over the DG solution. However, the errors are not as reduced as for the $3\pi/4$ orientation. The contour line plots in Figure 6 of the pointwise error profile show that there is a clear error reduction compared to the unfiltered solution.

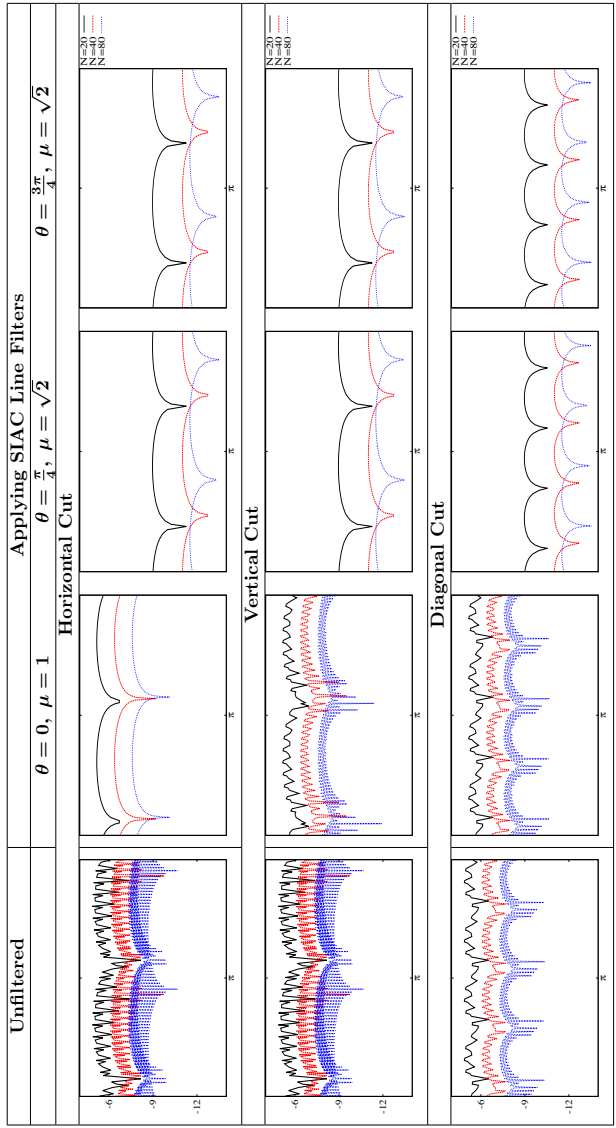

Figure 4: Pointwise error profile slices along the $x$-axis (horziontal), $y$-axis (vertical) and the diagonal directions when applying three different Line Filters on the DG solution to Problem (33) with initial condition $u_0(x,y) = \sin(x+y)$ and using $\mathbb{P}^3$ polynomials.

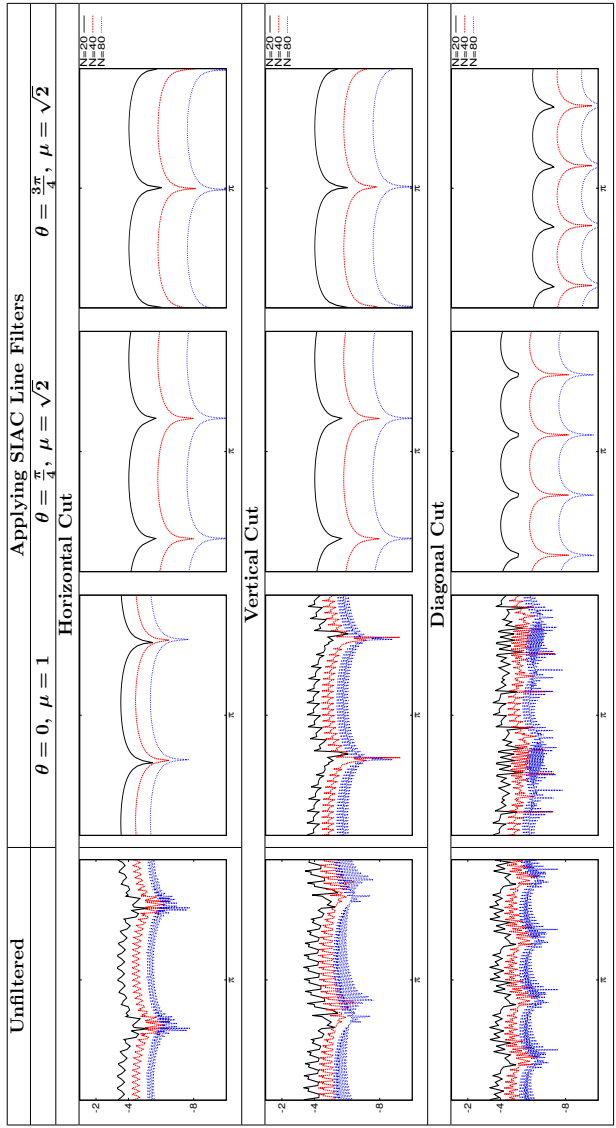

Figure 5: Point-wise error profile slices along the $x$-axis (horizontal), $y$-axis (vertical) and the diagonal directions when applying three different Line Filters on the DG solution to Problem (33) with initial condition $u_0(x, y) = \sin x \cos y$ and using $\mathbb{P}^2$ polynomials.

We conclude this numerical section studying the rotations for line filters alone to show the importance of the orientation for superconvergence extraction. For this, we considered the advection equation with initial condition $u(x,y) = \sin x \cdot \cos y$ and applied three line filters with rotations $\theta = 0, \pi/4$ and $3\pi/4$. Figure **??** shows the contour line error plots of the DG solution before and after filtering with the three line filters. Furthermore, Table 2 shows the $L^2$ errors and orders before and after line filtering. Observe that for the zero rotation, the filter is unable to raise the convergence rate, and remains the same as for the unfiltered solution. This is because the filtering is only done along the $x-$axis and therefore the derivatives in the $y-$direction are not included which invalidates the theory. On the other hand, we see that for rotation angles of $\theta \neq n\pi, n\pi/2, n \in \mathbb{N}$ the convergence rate is indeed raised. It is important to note that although the unrotated filter stays at $k+1$ accuracy, the global error is lower than the one for the unfiltered solution. However, this difference is not as remarkable as the error reduction for the rotated filter. For the choice of $\theta = \pi/4$, $3\pi/4$ we see clear error reduction, which is even more considerable when the mesh is refined or the polynomial degree of the approximation is increased. Furthermore, note that in Figure **??**, the rotation $\theta = 3\pi/4$ produces different error contours. We speculate that this rotation empathises the zeros of the sine function.

## 5 A note on computation

Finally, we would like to turn our attention to the computational challenges arising from the implementation of the SIAC filters. We will discuss the main filtering operations in terms of computational complexity and simulation times, highlighting the many advantages that line filtering has over the traditional 2D filtering.

Although the kernel has compact support and the convolution reduces to a small domain region, this integral contains several discontinuities. Hence, the total integral is computed as a sum of integrable regions. One of the most challenging and intense computations is to actually find and sort all these regions, delimited by the mesh element boundaries and the kernel breaks; the natural discontinuous structure of DG produces a solution that is integrable only inside the elements. Furthermore, since the kernel is built as a linear combination of B-Splines, there is only $k-1$ smoothness for each spline. More details on the computational aspects of the convolution and the impact on the error can be found in [13]. The integral of the filtering convolution was calculated using Gauss quadratures. For the 2D Cartesian axis aligned filter, the integral along each of the quadrilateral regions was computed as a tensor product of univariate quadratures along each direction. This technique is exact for polynomial integrands provided sufficient quadrature points are used. In Figure 8 we show total number of integral regions required to post-process one point using the tensor product filter aligned with the Cartesian axis and the $3\pi/4$-line filter. Table 3 and Table 4 summarize the number of operations and elapsed times required to post-process a single point using the same computer for both filters.

Notice that, both computational times and costs are notoriously reduced using a line filter. Line filters use one dimensional quadrature rules and the total number of integrals and quadrature sums match. On the other hand, the total number of sums using 2D filters is increased by a factor of $n^2$, where $n$ denotes the total number of integrals. The computational times shown in Table 4 not only show clear benefits for the line filter but also indicate a slower time increase as we raise the number and order of B-Splines. This represents a great advantage, as one important limiting factor on the applications of 2D SIAC filters is the long computational times of higher degree filters.

Table 1: $L^2$ errors and convergence rates before and after filtering the DG solution to Problem (33) with $u_0(x,y) = \sin(x+y)$ using the original 2D filter and two line filters. $h$ denotes the DG mesh size.

| | Unfiltered | | 2D Filtering $\theta=0, H=h$ | | Line Filtering $H = \mu h$ | $\theta = \pi/4$ | | $\theta = 3\pi/4$ | |
|---|---|---|---|---|---|---|---|---|---|
| N | $L^2$-Error | Order | $L^2$-Error | Order | $\mu$ | $L^2$-Error | Order | $L^2$-Error | Order |
| | | | | | $\mathbb{P}^1$ | | | | |
| 20 | 9.7e-03 | - | 1.6e-03 | - | 1 | 2.3e-03 | - | 2.e-03 | - |
| | | | | | $\sqrt{2}$ | 2.7e-03 | - | 1.5e-03 | - |
| 40 | 2.4e-03 | 2.02 | 2.0e-04 | 3.05 | 1 | 3.7e-04 | 2.62 | 3.6e-04 | 2.42 |
| | | | | | $\sqrt{2}$ | 2.6e-04 | 3.33 | 1.9e-04 | 2.98 |
| 80 | 5.9e-04 | 2.01 | - | - | 1 | 7.9e-05 | 2.22 | 8.1e-05 | 2.17 |
| | | | | | $\sqrt{2}$ | 2.8e-05 | 3.21 | 2.4e-05 | 2.99 |
| | | | | | $\mathbb{P}^2$ | | | | |
| 20 | 2.4e-04 | - | 6.1e-06 | - | 1 | 2.8e-05 | - | 1.5e-05 | - |
| | | | | | $\sqrt{2}$ | 1.4e-04 | - | 1.5e-06 | - |
| 40 | 2.9e-05 | 3.01 | 1.2e-07 | 5.71 | 1 | 2.0e-06 | 3.83 | 1.8e-06 | 3.02 |
| | | | | | $\sqrt{2}$ | 2.3e-06 | 5.91 | 4.7e-08 | 4.99 |
| 80 | 1.5e-05 | 3.01 | - | - | 1 | 2.2e-07 | 3.34 | 2.3e-07 | 3.01 |
| | | | | | $\sqrt{2}$ | 3.7e-08 | 5.95 | 1.5e-09 | 5.00 |
| | | | | | $\mathbb{P}^3$ | | | | |
| 20 | 4.5e-06 | - | 1.4-e-07 | - | 1 | 1.1e-06 | - | 7.1e-08 | - |
| | | | | | $\sqrt{2}$ | 1.6e-05 | - | 7.7e-10 | - |
| 40 | 8.3e-07 | 4.01 | 5.6e-10 | 7.96 | 1 | 6.8e-09 | 7.31 | 4.4e-09 | 4.02 |
| | | | | | $\sqrt{2}$ | 6.9e-08 | 7.87 | 6.9e-12 | 6.79 |
| 80 | 5.2e-08 | 4.00 | - | - | 1 | 2.8e-10 | 6.00 | 7.8e-11 | 5.81 |
| | | | | | $\sqrt{2}$ | 2.7e-10 | 7.97 | 2.9e-14 | 7.90 |

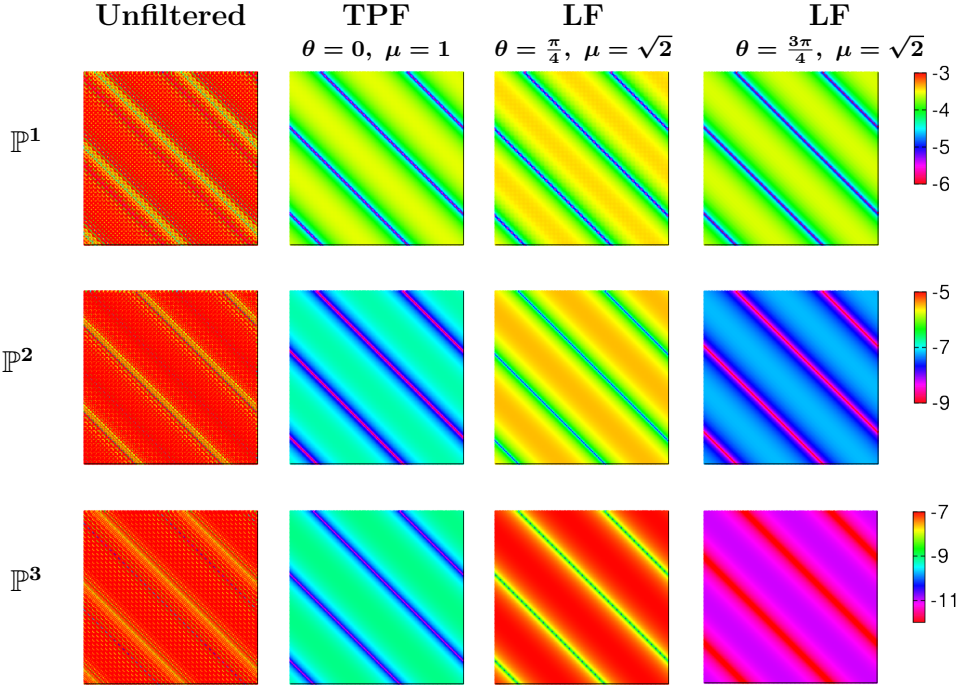

Figure 6: Contour line error plots (log) before and after filtering the DG solution to Problem (33) with $u_0(x,y) = \sin(x+y)$ using a Tensor Product Filter (TPF) and two Line Filters (LFs).

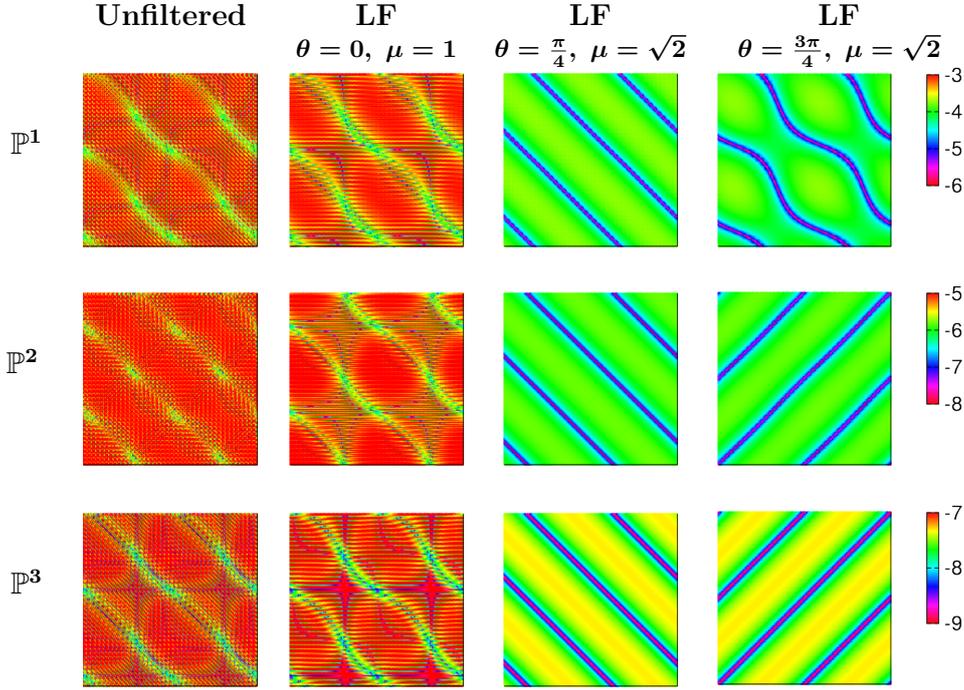

Figure 7: Contour line error plots (log) before and after filtering the DG solution to Problem (33) with $u_0(x,y) = \sin x \cdot \cos y$ applying three different Line Filters (LFs).

Table 2: $L^2$ errors and convergence rates before and after filtering the DG solution to Problem (33) with $u_0(x,y) = \sin x \cdot \cos y$ using three different line filters. $h$ denotes the mesh size.

| | Unfiltered | | Line Filtering | | | | | |
|---|---|---|---|---|---|---|---|---|
| | | | $H = \mu h$ | $\theta = 0$ | | $\theta = \pi/4$ | | $\theta = 3\pi/4$ | |
| N | $L^2$-Error | Order | $\mu$ | $L^2$-Error | Order | $L^2$-Error | Order | $L^2$-Error | Order |
| $\mathbb{P}^1$ | | | | | | | | | |
| 20 | 5.2e-03 | - | 1 | 3.7e-03 | - | 1.2e-03 | - | 1.0e-03 | - |
| | | | $\sqrt{2}$ | 3.7e-03 | - | 1.3e-03 | - | 9.7e-04 | - |
| 40 | 1.3e-03 | 2.02 | 1 | 9.1e-04 | 2.03 | 2.0e-04 | 2.57 | 2.0e-04 | 2.41 |
| | | | $\sqrt{2}$ | 9.1e-04 | 2.04 | 1.3e-04 | 3.33 | 1.0e-04 | 3.23 |
| 80 | 3.2e-04 | 2.01 | 1 | 2.2e-04 | 2.01 | 4.3e-05 | 2.19 | 4.4e-05 | 2.15 |
| | | | $\sqrt{2}$ | 2.2e-04 | 2.01 | 1.4e-05 | 3.21 | 1.2e-05 | 3.08 |
| $\mathbb{P}^2$ | | | | | | | | | |
| 20 | 1.3e-04 | - | 1 | 9.0e-05 | - | 1.4-e05 | - | 1.2e-05 | - |
| | | | $\sqrt{2}$ | 9.1e-05 | - | 6.8e-05 | - | 6.7e-05 | - |
| 40 | 1.6e-05 | 3.01 | 1 | 1.1e-05 | 3.01 | 1.0e-06 | 3.80 | 9.8e-07 | 3.61 |
| | | | $\sqrt{2}$ | 1.1e-05 | 3.02 | 1.1e-06 | 5.91 | 1.1e-06 | 5.92 |
| 80 | 2.0e-06 | 3.00 | 1 | 1.4e-06 | 3.00 | 1.2e-07 | 3.12 | 1.2e-07 | 3.03 |
| | | | $\sqrt{2}$ | 1.4e-06 | 3.00 | 1.8e-08 | 5.95 | 1.8e-08 | 5.98 |
| $\mathbb{P}^3$ | | | | | | | | | |
| 20 | 2.4e-06 | - | 1 | 1.7e-06 | - | 5.4e-07 | - | 5.4e-07 | - |
| | | | $\sqrt{2}$ | 1.9e-06 | - | 8.1-e06 | - | 8.1e-06 | - |
| 40 | 1.5e-07 | 4.01 | 1 | 1.1e-07 | 4.01 | 3.5e-09 | 7.28 | 3.2e-09 | 7.40 |
| | | | $\sqrt{2}$ | 1.1e-07 | 4.13 | 3.4e-08 | 7.87 | 3.4e-08 | 7.87 |
| 80 | 9.5e-09 | 4.00 | 1 | 6.7e-09 | 4.00 | 1.4e-10 | 4.59 | 1.4e-10 | 4.48 |
| | | | $\sqrt{2}$ | 6.7e-09 | 4.00 | 1.4e-10 | 7.97 | 1.4e-10 | 7.97 |

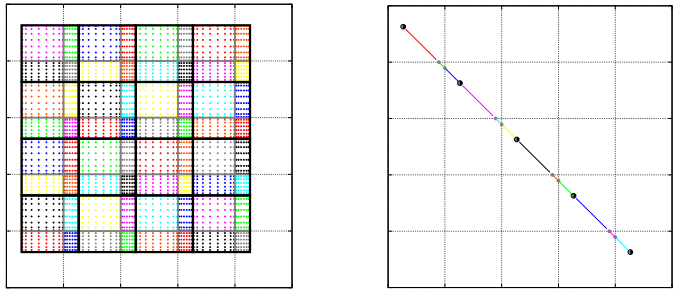

Figure 8: Total integrable subregions when post-processing a particular point with a Tensor Product Filter aligned with the Cartesian axis (left) and a rotated Line filter (right) over an uniform squared mesh.

| Filter Type | Intersection Scans | Integrals | Quadrature Sums |
|---|---|---|---|
| 2D No Rotation | 64 | 63 | 3969 |
| Line Filter | 4 | 12 | 12 |

Table 3: Summary of the number of operations required to compute the filtering convolution corresponding to the filters from Figure 8.

# 6 Conclusions and Future Work

In this article, we have introduced the SIAC Line Filter as a means of post-processing multi-dimensional data using a one-dimensional support. Introducing a rotation allows a more robust definition of the post-processor whilst conserving the properties defining these filters. By rotating the filter support, we have a better understanding of the relation between the filter, the underlying mesh and numerical scheme. The lower dimensional filtering technique has promising applications in fluid flow visualization of DG solutions. It has all the advantages of the original 2D Cartesian aligned filter while being more computationally efficient. We have developed theoretical error estimates showing how we can achieve similar results to the 2D Cartesian coordinate aligned SIAC filter. Numerical investigation demonstrated that the SIAC line filters provides more error reduction capability. The low computational costs associated to these filters makes them a very attractive tool for the scientific community. In the future we hope to apply line filters to flow visualization such as during streamline computations.

| No. of Splines and degree | Tensor Product Filter | Line Filter |
|---|---|---|
| 3, 1 | 0.42 | 0.03 |
| 5, 2 | 1.48 | 0.06 |
| 7, 3 | 3.8 | 0.1 |

Table 4: Computational times (seconds) taken by each filter from Figure 8 to post-process the same point.

# 7 Acknowledgements

The numerical data to which the filters were applied was generated through the Nektar++ software [3]. The first author would like to give special thanks to Chris Cantwell and David Moxey for all their help when adapting the Nektar++ code as well to Mike Kirby and the SCI Institute (University of Utah) for providing access to the school supercomputers during her visit.

# A Completing the proof in Lemma 2

We will show that the formula:

$$\partial_{u_\theta}^\alpha f(x,y) = \sum_{m=0}^{\alpha} \binom{\alpha}{m} \partial_{u_x^\theta}^{\alpha-m} \partial_{u_y^\theta}^m f(x - \frac{m}{2}\cos\theta, y + \frac{\alpha-m}{2}\sin\theta)$$

holds for any smooth function $f$ and $\alpha > 1$. The first order divided difference ($\alpha = 1$) was proven in Lemma 2.

Assume now that the formula holds for:

$$\partial_{u_\theta}^{\alpha-1} f(x,y) = \sum_{m=0}^{\alpha-1} \binom{\alpha-1}{m} \partial_{u_x^\theta}^{\alpha-1-m} \partial_{u_y^\theta}^m f(x - \frac{m}{2}\cos\theta, y + \frac{\alpha-1-m}{2}\sin\theta).$$

Then,

$$\partial_{u_\theta}^\alpha f(x,y) = \partial_{u_\theta} \left( \partial_{u_\theta}^{\alpha-1} f(x,y) \right)$$

$$= \sum_{m=0}^{\alpha-1} \binom{\alpha-1}{m} \partial_{u_\theta} \left( \partial_{u_x^\theta}^{\alpha-1-m} \partial_{u_y^\theta}^m f(x - \frac{m}{2}\cos\theta, y + \frac{\alpha-1-m}{2}\sin\theta) \right).$$

We know that

$$\partial_{u_\theta} f(x,y) = \partial_{u_x^\theta} f(x, y + \frac{1}{2}\sin\theta) + \partial_{u_y^\theta} f(x - \frac{1}{2}\cos\theta, y).$$

Then,

$$\partial_{u_\theta}^\alpha f(x,y) = (a) + (b),$$

$$(a) = \sum_{m=0}^{\alpha-1} \binom{\alpha-1}{m} \partial_{u_x^\theta}^{\alpha-m} \partial_{u_y^\theta}^m f(x - \frac{m}{2}\cos\theta, y + \frac{\alpha-m}{2}\sin\theta)$$

$$(b) = \sum_{m=0}^{\alpha-1} \binom{\alpha-1}{m} \partial_{u_x^\theta}^{\alpha-(m+1)} \partial_{u_y^\theta}^{m+1} f(x - \frac{m+1}{2}\cos\theta, y + \frac{\alpha-(m+1)}{2}\sin\theta).$$

Write the first term as:

$$(a) = \binom{\alpha-1}{0} \partial_{u_x^\theta}^\alpha f(x, y + \frac{\alpha}{2}\sin\theta)$$

$$+ \sum_{m=1}^{\alpha-1} \binom{\alpha-1}{m} \partial_{u_x^\theta}^{\alpha-m} \partial_{u_y^\theta}^m f(x - \frac{m}{2}\cos\theta, y + \frac{\alpha-m}{2}\sin\theta).$$

For the term $(b)$, we change $m \to m+1$,

$$(b) = \sum_{m=1}^{\alpha-1} \binom{\alpha-1}{m-1} \partial_{u_x^\theta}^{\alpha-m} \partial_{u_y^\theta}^m f(x - \frac{m}{2}H\cos\theta, y + \frac{\alpha-m}{2}H\sin\theta)$$

$$+ \binom{\alpha-1}{\alpha-1} \partial_{u_y^\theta}^\alpha f\left(x - \frac{\alpha}{2}\cos\theta, y\right).$$

Putting these two term together gives:

$$(a) + (b) = \sum_{m=1}^{\alpha-1} \underbrace{\binom{\alpha-1}{m} + \binom{\alpha-1}{m-1}}_{=\binom{\alpha}{m}} \partial_{u_x^\theta}^{\alpha-m} \partial_{u_y^\theta}^m f(x - \frac{m}{2}\cos\theta, y + \frac{\alpha-m}{2}\sin\theta)$$

$$+ \underbrace{\binom{\alpha-1}{0}}_{=\binom{\alpha}{0}} \partial_{u_x^\theta}^\alpha f(x, y + \frac{\alpha}{2}H\sin\theta) + \underbrace{\binom{\alpha-1}{\alpha-1}}_{=\binom{\alpha}{\alpha}} \partial_{u_y^\theta}^\alpha (x - \frac{\alpha}{2}H\cos\theta, y),$$

which gives the formula

$$\partial_{u_\theta}^\alpha f(x, y) = \sum_{m=0}^{\alpha} \binom{\alpha}{m} \partial_{u_x^\theta}^{\alpha-m} \partial_{u_y^\theta}^m f\left(x - \frac{m}{2}\cos\theta, y + \frac{\alpha-m}{2}\sin\theta\right).$$